\documentclass{article}
\usepackage{amsmath,amsfonts,amsthm,amssymb}
\DeclareMathSymbol\nullset{\mathord}{AMSb}{"3F}

\begin{document}
\newtheorem{theorem}{Theorem}[section]
\newtheorem{corollary}{Corollary}[section]
\newtheorem{lemma}{Lemma}[section]
\newtheorem{proposition}{Proposition}[section]
{\newtheorem{definition}{Definition}[section]}

\title{\centerline{\bf On a problem of A. V. Grishin}}
\author{C. Bekh-Ochir and S. A. Rankin
  }

\maketitle
\newcounter{parts}
\def\com#1,#2{[\,{#1},{#2}\,]}
\def\kx{k\langle X\rangle}
\def\kzerox{k_0\langle X\rangle}
\def\konex{k_1\langle X\rangle}
\def\set#1\endset{\{\,#1\,\}}
\def\rest#1{\,\hbox{\vrule height 6pt depth 9pt width .5pt}_{\,#1}}         
\def\choice#1,#2{\binom{#1}{#2}}
\def\kzerox{k_0\langle X\rangle}
\let\cong=\equiv
\def\comp{\mkern2mu\mathchoice%
        {\raise.35ex\hbox{$\scriptscriptstyle\circ$}}
        {\raise.35ex\hbox{$\scriptscriptstyle\circ$}}
        {\raise.14ex\hbox{$\scriptscriptstyle\circ$}}
        {\raise.14ex\hbox{$\scriptscriptstyle\circ$}}}
\def\from{\mkern2mu\hbox{\rm :}\mkern2mu}

\begin{abstract}
 In this note, we offer a short proof of V. V. Shchigolev's result that over any field $k$ of characteristic
 $p>2$, the $T$-space generated by $x_1^p,x_1^px_2^p\ldots$ is finitely based, which answered a question raised
 by A. V. Grishin. More precisely, we prove that for any field of any positive characteristic, $R_2^{(d)}=R_3^{(d)}$ 
 for every positive integer $d$, and that over an infinite field of characteristic $p>2$, $L_2=L_3$. Moreover,
 if the characteristic of $k$ does not divide $d$, we prove that $R_1^{(d)}$ is an ideal of $\kzerox$ and thus in
 particular, $R_1^{(d)}=R_2^{(d)}$. Finally, we show that 
 over any field of characteristic $p>2$, $R_1^{(d)}\ne R_2^{(d)}$ and $L_1\ne L_2$.
\end{abstract}

\section{Introduction}
 In \cite{Gr} (and later in \cite{ShGr}, the survery paper with V. V. Shchigolev), A. V. Grishin proved that in 
 the free associative algebra with countably infinite generating set $\set x_1,x_2,\ldots\endset$ over a field 
 of characteristic 2, the $T$-space that is generated by the set $\set x_1^2,x_1^2x_2^2,\ldots\endset$ 
 is not finitely based, and he raised the question as to 
 whether or not, in the corresponding setting but over a field
 of characteristic $p>2$, the $T$-space generatd by $\set x_1^p,x_1^px_2^p,\ldots\endset$ is finitely based.
 This was resolved by V. V. Shchigolev in \cite{Sh}, wherein he proved that over an infinite field
 of characteristic $p>2$, this $T$-space is finitely based. In fact, if we let $L_1$ denote the $T$-space
 generated by $\set x_1^p\endset$, and then for each positive integer $n$, let $L_{n+1}$ denote the $T$-space 
 generated by $L_n\cup L_n x_{n+1}^p$, Shchigolev proves in \cite{Sh} that $L_p=L_{p+1}$.
 To do this, he made use of another family of $T$-spaces defined in \cite{Sh} as follows. Let $k$ denote an
 arbitrary field of characteristic $p$, let $X=\set x_1,x_2,\ldots\endset$ be a countably infinite set, and 
 let $\kzerox$ denote the free associative $k$-algebra over the set $X$. For each positive integer
 $d$, let $S_d(x)$ denote the sum $\sum_{\sigma\in \Sigma_d} \prod_{i=1}^d x_{\sigma(i)}$, where $\Sigma_d$ is the
 symmetric group on $d$ letters. Let $R_1^{(d)}$ denote the $T$-space of $\kzerox$ that is generated by
 $S_d(x)$, and for each positive integer $n$, let $R_{n+1}^{(d)}$ denote the $T$-space of $\kzerox$ that is generated
 by $R_n^{(d)}\cup R_n^{(d)}S_d(x)$. As a key step in his demonstration that $L_p=L_{p+1}$, Shchigolev
 proves that for any positive integer $d$, $R_d^{(d)}=R_{d+1}^{(d)}$. This struck us as a bit curious -- why did the sequence
 $R_1^{(d)}\subseteq R_2^{(d)}\cdots R_d^{(d)}\subseteq R_{d+1}^{(d)}\subseteq \cdots$
 stabilize at the $d^{th}$ step? There did not seem to be a natural connection between the number of variables
 and the number of factors, and this led us to examine his argument more closely. The results
 of the original paper appear again in the survey paper \cite{ShGr} with some minor typographical errors corrected and in
 some cases, required conditions were clarified, but we note that there appears to be a minor error in the statement 
 of Lemma 15 of \cite{Sh} that did not get corrected in the survey paper.  
 Fortunately, this error does not affect the validity of the proof that $R_d^{(d)}=R_{d+1}^{(d)}$. Lemma 15 of
 \cite{Sh} states that for each $k=1,2,\ldots, d-1$, a certain polynomial $f_k$ is congruent modulo $R_1^{(d)}$
 to a summation expression. In fact, since it is not known whether or not $R_1^{(d)}$ is an ideal of 
 $\kzerox$, the best that can be said is that $f_k$ is congruent modulo $R_{k-1}^{(d)}$ to the summation expression.
 Shchigolev's proof that $R_d^{(d)}=R_{d+1}^{(d)}$ only needs that the summation expression
 that appears in the $d^{th}$ iterate belongs to $R_d^{(d)}$, and since the summation expression is congruent 
 modulo $R_{d-1}^{(d)}\subseteq R_d^{(d)}$ to $f_d$, and (it is apparent from the definition of
 $f_k$) $f_k\in R_k^{(d)}$ for each $k$, the desired conclusion holds.
 
 In this note, we offer a short proof that over any field $k$ (of any positive characteristic), $R_2^{(d)}=R_3^{(d)}$ 
 for every positive integer $d$, and that over an infinite field of characteristic $p>2$, $L_2=L_3$. Moreover,
 if the characteristic of $k$ does not divide $d$, we prove that $R_1^{(d)}$ is an ideal of $\kzerox$ and thus in
 particular, $R_1^{(d)}=R_2^{(d)}$. Finally, we prove that for every prime $p>2$, and any field $k$ of characteristic 
 $p$, $R_1^{(p)}\ne R_2^{(p)}$ and that $L_1\ne L_2$ (we remark that for an infinite field $k$, Shchigolev's argument 
 in \cite{Sh} can be used to imply that $R_1^{(p)}\ne R_2^{(p)}$).

 \section{$R_2^{(d)}=R_3^{(d)}$}
  For this section, $k$ is an arbitrary field. The proof of the first result is immediate.
  
  \begin{lemma}\label{lemma: basic}
    Let $d$ be a positive integer. Then 
    \begin{align*}
     S_{d+1}(x)&=\sum_{i=1}^{d+1}S_d(x_1,x_2,\ldots,\hat{x}_i,\ldots,x_{d+1})x_i\tag{1}\\ 
     &=S_d(x_1,x_2,\ldots,x_d)x_{d+1}+\sum_{i=1}^d S_d(x_1,x_2,\ldots,x_{d+1}x_i,\ldots,x_d)\tag{2}\\
     &=x_{d+1}S_d(x_1,x_2,\ldots,x_d)+\sum_{i=1}^d S_d(x_1,x_2,\ldots,x_ix_{d+1},\ldots,x_d).\tag{3}
    \end{align*}
  \end{lemma}
  
  \begin{corollary}\label{corollary: r1 commutator}
   For any $u\in A$, and any positive integer $d$, $\com S_d(x),u\in R_1^{(d)}$.
  \end{corollary}
  
  \proof
   This follows directly from (2) and (3) of Lemma \ref{lemma: basic}.
  \endproof
  
  We remark that in \cite{Sh}, Shchigolev proves that if the field is infinite, then for any $T$-space $L$, 
  if $v\in L$, then $\com v,u\in L$ for any $u\in A$.
  
  \begin{corollary}\label{corollary: mod r1}
   For any positive integer $d$, $S_{d+1}(x)\cong S_d(x)x_{d+1}\cong x_{d+1}S_d(x)\mod{R_1^{(d)}}$.
  \end{corollary}
  
  \proof
   This is also immediate from (2) and (3) of Lemma \ref{lemma: basic}.
  \endproof
  
  \begin{corollary}\label{corollary: dsd}
   For any positive integer $d$, $d\,S_{d+1}(x)\in R_1^{(d)}$.
  \end{corollary}
  
  \proof
   Note that $S_d(x_1,x_2,\ldots,x_d)= S_d(x_{\sigma(1)},x_{\sigma(2)},\ldots,x_{\sigma(d)})$
   for any positive integer $d$ and any $\sigma\in \Sigma_d$. By Corollary \ref{corollary: mod r1}, applied $d+1$
   times with a different variable pulled out each time, we obtain 
   $$
    (d+1)S_{d+1}(x)\cong\sum_{i=1}^{d+1}S_d(x_1,x_2,\ldots,\hat{x}_i,\ldots,x_{d+1})x_i\mod{R_1^{(d)}},
   $$
   and so the result follows from (1) of Lemma \ref{lemma: basic}.
  \endproof
  
  If the characteristic of $k$ does not divide $d$, it follows from Corollary \ref{corollary: dsd} that 
  $S_{d+1}(x)\in R_1^{(d)}$, and then Corollary \ref{corollary: mod r1} implies that $R_1^{(d)}$ is an ideal of $\kzerox$.
  In particular, if the characteristic of $k$ does not divide $d$, then $R_1^{(d)}=R_2^{(d)}$.
  
 \begin{proposition}\label{proposition: r2=r3}
  Let $k$ be a field of characteristic $p$, and let $d$ be a positive multiple of $p$. Then
  $R_2^{(d)}=R_3^{(d)}$.
 \end{proposition}

 \proof
  By (1) of Lemma \ref{lemma: basic}, $S_d(x)x_{d+1}+
  \sum_{i=1}^{d}S_d(x_1,x_2,\ldots,\hat{x}_i,\ldots,x_{d+1})x_i = S_{d+1}(x)$, and by
  Corollary \ref{corollary: mod r1}, we have $\sum_{i=1}^{d}S_d(x_1,x_2,\ldots,\hat{x}_i,\ldots,x_{d+1})x_i
  =S_{d+1}(x)-S_d(x)x_{d+1}\in R_1^{(d)}$. Let $u\in R_1^{(d)}$. Then 
  \begin{align*}
  uS_d(x_2,\ldots,x_{d+1})x_1 + &\sum_{i=2}^{d} uS_d(x_1,\ldots,\hat{x}_i,\ldots,x_{d+1})x_i\\
  &\hskip4pt=u\sum_{i=1}^{d}S_d(x_1,x_2,\ldots,\hat{x}_i,\ldots,x_{d+1})x_i\in R_2^{(d)}\mkern-5mu\hbox to 0pt{.\hss}
  \end{align*}
  Since $u\in R_1^{(d)}$, it follows from Corollary \ref{corollary: r1 commutator} that
  for each $i=2,\ldots,d$, 
  $$
   uS_d(x_1,\ldots,\hat{x}_i,\ldots,x_{d+1})x_i \cong 
  S_d(x_1,\ldots,\hat{x}_i,\ldots,x_{d+1})x_iu \mod{R_1^{(d)}},
  $$
  and by two applications of Corollary \ref{corollary: mod r1}, we then obtain
  \begin{align*}
  S_d(x_1,\ldots,\hat{x}_i,\ldots,x_{d+1})x_iu &\cong S_{d+1}(x_1,\ldots,\hat{x}_i,\ldots,x_{d+1},x_iu)\\
  &\cong S_d(x_2,\ldots,\hat{x}_i,\ldots,x_{d+1},x_iu)x_1\mod{R_1^{(d)}}.
  \end{align*}
  Thus for any $v\in R_1^{(d)}$, upon replacing $x_1$ by $v$ we obtain that 
  $$
  uS_d(x_2,\ldots,x_{d+1})v +\bigl(\sum_{i=2}^{d} S_d(x_2,\ldots,\hat{x}_i,\ldots,x_{d+1},x_iu)\bigr)v\in R_2^{(d},
  $$
  and so $uS_d(x_2,\ldots,x_{d+1})v\in R_2^{(d)}$ for all $u,v\in R_1^{(d)}$. It follows that $R_3^{(d)}\subseteq R_2^{(d)}$.
 \endproof 
   
\section{$L_2=L_3$}
   
The central idea behind Shchigolev's proof that $L_p=L_{p+1}$ is encapsulated in the following lemma.

\begin{lemma}\label{lemma: fundamental}
 Let $k$ be an infinite field of characteristic $p>2$, $i$ be any positive integer, and $u$ a multihomogeneous 
 element of $L_i$. 
 If $uS_p(x)S_p(y)\in L_{i+1}$, where $u$, $S_p(x)$, and $S_p(y)$ have no generators of $\kzerox$ in common,
 then $uz_1^pz_2^p\in L_{i+1}$, where $z_1\ne z_2$ and neither appears in $u$.
\end{lemma}

\proof
 For convenience, for any positive integer $n$, and any subset $U$ of $J_n=\set 1,2,\ldots,n\endset$,
 let 
 $$
   X_n(U)=\prod_{i=1}^n x_i\rest{x_i=\hbox{$\scriptstyle\begin{cases} z_1 & i\in U\\ \vspace{-4pt} z_2 & i\notin U\end{cases}$}}.
 $$  
 Then $uS_p(x)S_p(y)\in L_{i+1}$ implies that for every $j$ with $1\le j\le p-1$, we have
 \begin{align*}
  u(j!(p-j)!)^2&\sum_{\substack{ U\subseteq J_p\\|U|=j}}\mkern-7mu X_p(U)\mkern-10mu\sum_{\substack{U\subseteq J_p\\|U|=p-j}} \mkern-10mu X_p(U)\\
  &\hskip10pt=uS_p(\underbrace{z_1,z_1,\ldots,z_1}_j,z_2,\ldots,z_2)S_p(\underbrace{z_1,z_1,\ldots,z_1}_{p-j},z_2,\ldots,z_2)
 \end{align*}
 so $u\sum_{\substack{U\subseteq J_p\\|U|=j}}\mkern-7mu X_p(U)\sum_{\substack{U\subseteq J_p\\|U|=p-j}}\mkern-10mu X_p(U)\in L_{i+1}$,
 and thus 
 $$
   g=\sum_{j=1}^{p-1}u\sum_{\substack{U\subseteq J_p\\|U|=j}} X_p(U)\sum_{\substack{U\subseteq J_p\\|U|=p-j}} X_p(U)
 $$
 is an element of $L_{i+1}$. On the other hand, since $u\in L_i$, we have $u(z^2)^p\in L_{i+1}$, and so
 $u(z_1+z_2)^{2p}\in L_{i+1}$. As $k$ is infinite, $h$, the sum of all multihomogeneous component of $u(z_1+z_2)^{2p}$ with
 degree $p$ for each of $z_1$ and $z_2$, belongs to $L_{i+1}$. We have
 $$
  h= u \sum_{\substack{ U\subseteq J_{2p}\\|U|=p}}\mkern-7mu X_{2p}(U)
  =u(z_1^pz_2^p+z_2^pz_1^p)+g\in L_{i+1},
 $$
 and thus $u(z_1^pz_2^p+z_2^pz_1^p)\in L_{i+1}$. Furthermore, since $k$ is infinite and $z_1^p\in L_1$, we have $\com z_1^p,{z_2^p}
 \in L_1$ and thus, since $u\in L_i$, $u(z_1^pz_2^p-z_2^pz_1^p)\in L_{i+1}$. It follows that $2uz_1^pz_2^p\in L_{i+1}$,
 and since $p>2$, we obtain $uz_1^pz_2^p\in L_{i+1}$, as required.
\endproof
 
 We need one additional fact.

\begin{lemma}\label{lemma: power}
 Let $k$ be an infinite field, and let $L$ be a $T$-space of $\kzerox$. For any positive integer $d$, and any
 $u\in L$, if $uz^d\in L$, where $z$ is a generator not appearing in $u$, then $uS_d(x)\in L$.
\end{lemma}

\proof
 Since $k$ is infinite, we may linearize $uz^d$ with respect to $z$ to obtain that $uS_d(x)\in L$.
\endproof

Lemma \ref{lemma: power} has the following corollary as an immediate consequence.

\begin{corollary}\label{corollary: rklk}
  Let $k$ be an infinite field of characteristic $p$. Then for every positive integer $k$, 
  $R_k^{(p)}\subseteq L_k$.
\end{corollary}

\begin{theorem}
 Let $k$ be an infinite field of characteristic $p>2$. Then $L_2=L_3$.
\end{theorem}

\proof
 By Corollary \ref{corollary: rklk}, $R_2^{(p)}\subseteq L_2$, and by Proposition \ref{proposition: r2=r3}, $R_3^{(p)}=R_2^{(p)}$.
 Thus $S_p(x)S_p(y)S_p(z)\in R_3=R_2\subseteq L_2$. As well, $S_p(x)\in R_1^{(p)}\subseteq L_1$,
 so by Lemma \ref{lemma: fundamental}, $S_p(x)z_1^pz_2^p\in L_2$. But now, since $S_p(x)z_1^p\in L_2$,
 we have by Lemma \ref{lemma: power} that $S_p(x)z_1^pS_p(y)\in L_2$. Since $S_p(x)\in R_1^{(d)}\subseteq L_1\subseteq
 L_2$, we have $\com S_p(x),{z_1^pS_p(y)}\in L_2$ and thus $z_1^pS_p(y)S_p(x)\in L_2$. As $z_1^p\in L_1$,
 we may apply Lemma \ref{lemma: fundamental} again to obtain that $z_1^pz_2^pz_3^p\in L_2$. Thus $L_3\subseteq L_2$.
\endproof   

\section{A study of $R_1^{p}$ for prime $p>2$}

In this section, we explore the structure of $R_1^{(p)}$ for an arbitrary prime $p$ and an arbitrary field of
characteristic $p$.

\begin{definition}
 For any positive integer $n$, and any $i$ with $1\le i\le n+1$, let $M_i^{n+1}=\prod_{j=1}^{n+1} a_j$,
 where $a_j=x$ if $j=i$, otherwise $a_j=y$.
\end{definition} 

The proof of the following result involves an elementary inductive argument based on Pascal's identity,
and has been omitted.

\begin{lemma}\label{lemma: basic comm identity}
 For any integer $n\ge1$, $[x,\underbrace{y,\ldots,y}_{n}] = \sum_{i=1}^{n+1} (-1)^{i+1}\choice n,{i-1} M_i^{n+1}$.
\end{lemma}

\begin{lemma}\label{lemma: p identity}
 For any prime $p$, and any integer $i$ with $0\le i\le p-1$, $\choice p-1,i\cong (-1)^i\mod{p}$.
\end{lemma}

\proof
 The result is immediate for $p=2$, while for $p>2$, it follows also from Pascal's identity $\choice p-1,i+\choice p-1,{i-1}=\choice p,i$,
 which is zero modulo $p$ if $0<i<p$. Thus $\choice p-1,i\cong -\choice p-1,{i-1}$ for $i=1,2,\ldots,p-1$. The result follows
 by induction based on $\choice p-1,{p-1}=1$.
\endproof 
 
\begin{corollary}\label{corollary: comm sum identity}
 Let $p$ be a prime. Then $[x,\underbrace{y,\ldots,y}_{p-1}]=\sum_{i=1}^p (-1)^{i+1}\choice p-1,{i-1} M_i^p$.
\end{corollary}

\proof
 By Lemma \ref{lemma: basic comm identity}, we have $[x,\underbrace{y,\ldots,y}_{p-1}] = 
 \sum_{i=1}^{p} (-1)^{i+1}\choice p-1,{i-1} M_i^{p}$, and thus by Lemma \ref{lemma: p identity},
 $[x,\underbrace{y,\ldots,y}_{p-1}] =  \sum_{i=1}^{p} (-1)^{i+1}(-1)^{i-1} M_i^{p}$,
 as required.
\endproof

\begin{proposition}\label{proposition: sp result}
 For any prime $p$, $S_p(x,\underbrace{y,\ldots,y}_{p-1}) = -[x,\underbrace{y,\ldots,y}_{p-1}]$.
\end{proposition} 

\proof
 $S_p(x,\underbrace{y,\ldots,y}_{p-1}) = (p-1)!\sum_{i=1}^p M_i^p$, and so by Wilson's theorem and
 Corollary \ref{corollary: comm sum identity}, we obtain the desired result.
\endproof

 For any integer $n\ge 2$, let $\Sigma_n^*$ denote the permutation group 
 on $\set 2,3,\ldots,n\endset$.
 
\begin{corollary}\label{corollary: fund sp identity}
 $S_p(x_1,x_2,\ldots,x_p)=\sum_{\sigma\in \Sigma_p^*}[x_1,x_{\sigma(2)},\ldots,x_{\sigma(p)}]$ 
 for any prime $p$.
\end{corollary}

\proof
 Let $u=\sum_{j=2}^p x_j$. Then $S_p(x_1,\underbrace{u,\ldots,u}_{p-1}) = -[x_1,\underbrace{u,\ldots,u}_{p-1}]$.
 Upon expansion, we obtain that $(p-1)!s_p(x_1,x_2,\ldots,x_p)+\sum_{\sigma\in \Sigma_p^*}[x_1,x_{\sigma(2)},\ldots,x_{\sigma(p)}]$,
 a sum of monomials each of which depends essentially on $x_1,x_2,\ldots,x_p$,
 is equal to a sum of monomials each of which is missing at least one of the variables $x_2,\ldots,x_p$.
 Since the set of all monomials forms a linear basis for $\kzerox$, it follows that
 $(p-1)!S_p(x_1,x_2,\ldots,x_p)+\sum_{\sigma\in \Sigma_p^*}[x_1,x_{\sigma(2)},\ldots,x_{\sigma(p)}]=0$. Since $(p-1)!=-1$,
 the result follows.
\endproof

\begin{proposition}\label{proposition: nice}
 $R_1^{(p)}=\set [x,\underbrace{y,\ldots,y}_{p-1}]\endset^S$ for any prime $p$ and any field of characteristic $p$,
\end{proposition}

\proof
 By Proposition \ref{proposition: sp result}, $\set [x,\underbrace{y,\ldots,y}_{p-1}]\endset^S\subseteq R_1^{(p)}$. On the other hand,
 if we let $u=\sum_{i=2}^p x_i$, then 
 $S_p(x_1,\underbrace{u,\ldots,u}_{p-1})=(p-1)!S_p(x_1,x_2,\ldots,x_p)+w$, where $w$
 is a sum of monomials each of which is missing at least one variable, and by Proposition \ref{proposition: sp result},
 $S_p(x_1,\underbrace{u,\ldots,u}_{p-1})\in \set [x,\underbrace{y,\ldots,y}_{p-1}]\endset^S$. Thus 
 $(p-1)!S_p(x_1,x_2,\ldots,x_p)$, the part of $S_p(x_1,\underbrace{u,\ldots,u}_{p-1})$ which depends essentially on
 $x_1,x_2,\ldots,x_p$ belongs to $\set [x,\underbrace{y,\ldots,y}_{p-1}]\endset^S$ as well, and so we have
 $R_1^{(p)}\subseteq \set [x,\underbrace{y,\ldots,y}_{p-1}]\endset^S$.
\endproof  

 \section{$R_1^{(p)}\ne R_2^{(p)}$ and $L_1\ne L_2$}

  \begin{lemma}\label{lemma: basic coefficient}
    For any $x,y\in X$ and monomials $u_1,u_2,\ldots,u_p$ in $\kzerox$, the coefficients of $(xy)^p$ and
    of $(yx)^p$ in $S_p(u_1,\ldots,u_p)$ sum to zero.
  \end{lemma}

  \begin{proof}
   We must show that if there exists $\sigma\in \Sigma_p$ such that $\prod_{i=1}^p u_{\sigma(i)}=(xy)^p$ or
   $(yx)^p$, then the number of such permutations in $\Sigma_p$ is a multiple of $p$ (actually, this argument
   does not depend on $p$ being prime). Let $\gamma\in\Sigma_p$ denote the cyclic permutation that sends $i$ to $i+1$
   for $1\le i\le p-1$, with $p$ being sent to 1. Then the equivalence relation defined on $\Sigma_p$ by
   saying $\sigma$ is related to $\tau$ if and only if there exists an integer $t$ such that $\sigma=\gamma^t\comp \tau$ 
   has as its equivalence classes the right cosets of the cyclic subgroup generated by $\gamma$, so each equivalence
   class has size $p$. Suppose that that $\sigma\in \Sigma_p$ is such that $\prod_{i=1}^p u_{\sigma(i)}=(xy)^p$ or
   $(yx)^p$. Then $\prod_{i=1}^p u_{\gamma\comp\sigma(i)}=(xy)^p$ or $(yx)^p$ as well (depending on whether
   $u_{\sigma(p)}$ starts with $x$ or $y$), and thus for every $\tau$ in the equivalence class of $\sigma$, we have
   $\prod_{i=1}^p u_{\tau(i)}=(xy)^p$ or $(yx)^p$. Thus the permutations of $u_1,\ldots,u_p$ that produce either
   $(xy)^p$ or $(yx)^p$ can be partitioned into cells of size $p$, and so the result follows.
  \end{proof}

\begin{theorem}\label{theorem: r1r2}
 If $p>2$, then $R_1^{p}\ne R_2^{(p)}$.
\end{theorem}

\begin{proof} 
 Suppose to the contrary that $R_1^{p}= R_2^{(p)}$. Then in particular, 
 $$
 w=S_p(\underbrace{x,x,\ldots,x}_{\frac{p+1}{2}},\underbrace{y,y,\ldots,y}_{\frac{p-1}{2}})
    S_p(\underbrace{y,y,\ldots,y}_{\frac{p+1}{2}},\underbrace{x,x,\ldots,x}_{\frac{p-1}{2}})\in R_1^{(p)}
 $$
 (note that this is where we use the fact that $p$ is odd).
 Since $w\in R_1^{(p)}$, $w$ can be written as a linear combination of terms of the form $S_p(u_1,u_2,\ldots,u_p)$,
 $u_1,u_2,\ldots,u_p\in \kzerox$. As $S_p$ is multilinear, it follows that $w$ can be written as a linear combination
 of terms of the form $S_p(u_1,u_2,\ldots,u_p)$, where $u_1,u_2,\ldots,u_p$ are monomials in $\kzerox$, and thus by
 Lemma \ref{lemma: basic coefficient}, the sum of the coefficient of $(xy)^p$ and $(yx)^p$ in $w$ is zero. However,
 the coefficient of $(xy)^p$ in $w$ is $\displaystyle \bigl((\frac{p+1}{2})!(\frac{p-1}{2})!\bigr)^2\not\cong 0\mod{p}$, and
 the coefficient of $(yx)^p$ in $w$ is 0, which means that the sum of the coefficients of $(xy)^p$ and of $(yx)^p$ in
 $w$ is not zero, a contradiction.
\end{proof} 

\begin{proposition}\label{proposition: r1 not l1}
 $R_1^{(p)}\ne L_1$.
\end{proposition}

\begin{proof}
 Suppose to the contrary that $R_1^{(p)}=L_1$. Then $(xy)^p\in R_1^{(p)}$, which means that $(xy)^p$ can be
 written as a linear combination of terms of the form (using the fact that $S_p$ is multilinear) $S_p(u_1,u_2,\ldots,u_p)$,
 where $u_1,u_2,\ldots,u_p$ are monomials in $\kzerox$. But by Lemma \ref{lemma: basic coefficient}, in any linear combination of
 such terms, the sum of the coefficients of $(xy)^p$ and of $(yx)^p$ is zero, while the sum must in fact be 1. Thus
 we have a contradiction, and the result follows.
\end{proof}

The following result is well-known.

\begin{lemma}\label{lemma: up}
 For any $u,v\in\kzerox$ and integer $i$ with $1\le i\le p-1$, let 
 $$
  S_p(u,v;i)= S_p(\underbrace{u,u,\ldots,u}_{i},\underbrace{v,v,\ldots,v}_{p-i}).
 $$
 Then $\displaystyle (u+v)^p=u^p+v^p+\sum_{i=1}^{p-1}\frac{1}{i!(p-i)!}S_p(u,v;i)$.
\end{lemma} 

\begin{theorem}\label{theorem: l1l2}
 For any prime $p$, $L_1\ne L_2$.
\end{theorem}

\begin{proof}
 A. V. Grishin proved in \cite{Gr} that over an infinite field of characteristic 2, the $T$-space generated
 by $\set x_1^p,x_1^px_2^p,\ldots\endset$ is not finitely based, and thus $L_i\ne L_j$ for any distinct $i$
 and $j$. We now consider the situation when $p>2$.
 Suppose to the contrary that $L_1=L_2$. Since $R_2^{(p)}\subseteq L_2$, we obtain that $R_2^{(p)}\subseteq L_1$ and
 so in particular, the element $w$ introduced in the proof of Theorem \ref{theorem: r1r2} belongs to $L_1$. But then
 $w$ is a linear combination of terms of the form $u^p$, $u\in \kzerox$. By Lemma \ref{lemma: up}, each term of the
 form $u^p$, $u\in \kzerox$, can be written as a linear combination of terms of the form $v^p$ or 
 $S_p(u_1,u_2,\ldots,u_p)$, where $v, u_1,u_2,\ldots,u_p$ are monomials in $\kzerox$. We note that $w$ is multihomogeneous
 of degree $p$ in each of $x$ and $y$,
 and as each expresssion of the form $S_p(u_1,u_2,\ldots,u_p)$ with $u_1,u_2,\ldots,u_p$ monomials in $\kzerox$ is
 multihomogenous, it follows that $w$ can be written as a multihomogeneous linear combination of terms of the form
 $(xy)^p$, $(yx)^p$, and $S_p(u_1,u_2,\ldots,u_p)$, where $u_1,u_2,\ldots,u_p$ are monomials in $\kzerox$.
 Thus we may assume that $w=r(xy)^p+s(yx)^p + z$, where $r,s\in k$ and $z$ is a linear combination
 of terms of the form $S_p(u_1,u_2,\ldots,u_p)$, with $u_1,u_2,\ldots,u_p$ monomials in $\kzerox$.
 Let $\alpha\from \kzerox\to \konex$ denote the algebra homomorphism that is determined by sending $x$ to itself, $y$ to 1, 
 and all other elements of $X$ to zero. Since the image of $\alpha$ is contained in the commutative ring $k<\set x\endset>$, 
 all elements of $\kzerox$ of the form $S_p(u_1,u_2,\ldots,u_p)$ will be mapped to zero. In particular, $\alpha(w)=0$
 and $\alpha(z)=0$, so we obtain that $(r+s)x^p=0$ and thus $r+s=0$. But this, together with Lemma \ref{lemma: basic coefficient},
 means that the sum of the coefficient of $(xy)^p$ and the coefficient of $(yx)^p$ in $w$ is equal to 0. As we have
 already observed in the proof of Theorem \ref{theorem: r1r2}, this is not the case, and so $L_1\ne L_2$.
\end{proof}

\end{document}